\documentclass{amsart}
\usepackage{amssymb,euscript,amsmath, mathrsfs}
\usepackage[dvips]{graphicx}
\usepackage[dvips]{color}

\newcounter{ENUM}
\newcommand{\itm}{\item}
\newenvironment{ilist}{\renewcommand{\theENUM}{\roman{ENUM}}\renewcommand{\itm}{\addtocounter{ENUM}{1}\item[(\theENUM)]}\begin{itemize}\setcounter{ENUM}{0}}{\end{itemize}}
\newenvironment{Ilist}{\renewcommand{\theENUM}{\Roman{ENUM}}\renewcommand{\itm}{\addtocounter{ENUM}{1}\item[(\theENUM)]}\begin{itemize}\setcounter{ENUM}{0}}{\end{itemize}}

\def\risom{\overset{\sim}{\rightarrow}}

\input xy
\xyoption{all}
\CompileMatrices

\def\ZZ{{\mathbb Z}}
\def\NN{{\mathbb N}}
\def\AA{{\mathbb A}}

\def\cE{{\mathscr E}}
\def\cF{{\mathscr F}}
\def\cG{{\mathscr G}}

\def\cO{{\mathscr O}}

\def\cLG{\mathcal{LG}}

\def\fp{{\mathfrak p}}

\def\ox{{\overline{x}}}
\def\talpha{{\tilde \alpha}}
\def\tbeta{{\tilde \beta}}

\def\vp{\varphi}

\def\Hom{\operatorname{Hom}}

\def\Spec{\operatorname{Spec}}

\def\im{\operatorname{im}}

\def\red{\operatorname{red}}

\def\depth{\operatorname{depth}}
\def\pr{\operatorname{pr}}

\newcommand{\colspan}{\operatorname{colspan}}
\newcommand{\Id}{\operatorname{Id}}

\newcommand{\margh}[1]{}

\newtheorem{thm}{Theorem}[section]

\newtheorem{lem}[thm]{Lemma}

\theoremstyle{definition}
\newtheorem{defn}[thm]{Definition}

\theoremstyle{remark}
\newtheorem{notn}[thm]{Notation}
\newtheorem{rem}[thm]{Remark}

\numberwithin{equation}{section}

\begin{document}
\title{Flatness of the linked Grassmannian}
\author{David Helm and Brian Osserman}
\begin{abstract} We show that the linked Grassmannian scheme, which arises
in a functorial compactification of spaces of limit linear series, and
in local models of certain Shimura varieties, is Cohen-Macaulay, reduced,
and flat. We give an application to spaces of limit linear series.
\end{abstract}
\thanks{The first author was supported by a fellowship from the NSF, and
the second author was supported by fellowships from the Clay Mathematics 
Institute and the NSF during the preparation of this paper.}
\maketitle

\section{Introduction}

The linked Grassmannian was introduced in \cite{os8} as a tool for 
the construction of a moduli scheme of limit linear series, providing a
functorial description which also compactified the spaces introduced by
Eisenbud and Harris in \cite{e-h1}. It also arises in other contexts:
for instance, the length $2$ case over $\ZZ_p$ gives a local model for 
certain unitary Shimura varieties \cite{go2}. 

A few basic properties of the linked Grassmannian were proved in \cite{os8},
but only as much as was necessary for the application to limit linear series.
In particular, although the dimensions were computed, a substantial question
that remained open (Question A.18 of {\it loc.\ cit.}) was the flatness of 
linked Grassmannian schemes over their 
base. In this paper, we address this question. We perform a detailed analysis
in the case of linked Grassmannians of length $2$, giving a local description
at each point (Theorem \ref{calc} below). This description allows us to 
invoke the work of Zhang \cite{zh1} to conclude that the linked 
Grassmannians are Cohen-Macaulay, from which we conclude flatness, as well
as reducedness. Finally, we use an inductive argument to conclude the same 
statements for linked Grassmannians of arbitrary length 
(Theorem \ref{lg-flat}).

As consequences, we are able to conclude (Theorem \ref{main-flat}) that when 
a limit linear series 
space has the expected dimension, it is flat over its base and 
Cohen-Macaulay. We also see that the linked Grassmannian itself
provides a new class of flat degenerations of the Grassmannian; such 
degenerations have frequently been of interest in the combinatorial study
of the Grassmannian, and the linked Grassmannian thus provides a potentially
useful new tool for such analysis. 

There is some substantial overlap of our 
results, in the special case of length $2$ over $\ZZ_p$, with reducedness
and flatness results of Goertz \cite{go2}.
However, because we are working over a more general base, and our local 
calculation is rather concise and gives a stronger Cohen-Macaulayness 
result, we have chosen to keep our presentation self-contained.

\section*{Acknowledgements} We would like to thank David Eisenbud, for 
explaining Lemma \ref{cm-int} to us.

\section{Review of the linked Grassmannian}

We briefly review the basic definitions and results of the linked
Grassmannian.

\begin{defn}\label{lgdef} Let $S$ be an integral, locally Cohen-Macaulay
scheme, and $\cE_1, \dots, \cE_n$ vector
bundles on $S$, each of rank $d$. Given maps
$f_i:\cE_i \rightarrow \cE_{i+1}$ and $g_i: \cE_{i+1} \rightarrow
\cE_i$, and a positive integer $r<d$, we denote by $\cLG:=\cLG(r,\{\cE_i\}_i,
\{f_i, g_i\}_i)$ the functor associating to each $S$-scheme $T$ the set of
sub-bundles $V_1, \dots, V_n$ of $\cE_{1,T}, \dots, \cE_{n,T}$ having rank
$r$ and
satisfying $f_{i,T}(V_i) \subset V_{i+1}$, $g_{i,T}(V_{i+1}) \subset V_i$
for all $i$.

We say that $\cLG$ is a {\bf linked Grassmannian} functor if the following
further conditions on the $f_i$ and $g_i$ are satisfied:

\begin{Ilist}
\itm There exists some $s \in \cO_S$ such that $f_i g_i = g_i f_i$ is
scalar multiplication by $s$ for all $i$.
\itm Wherever $s$ vanishes, the kernel of $f_i$ is precisely equal to
the image of $g_i$, and vice versa. More precisely, for any $i$ and given any
two integers $r_1$ and $r_2$ such that $r_1 + r_2 < d$, then the closed
subscheme of $S$ obtained as the locus where $f_i$ has rank less than or
equal to $r_1$ and $g_i$ has rank less than or equal to $r_2$ is empty.
\itm At any point of $S$, $\im f_i \cap \ker f_{i+1}=0$, and $\im
g_{i+1} \cap \ker g_i = 0$. More precisely, for any integer $r_1$, and any
$i$, we have locally closed subschemes of $S$ corresponding to the locus
where $f_i$ has rank exactly $r_1$, and $f_{i+1} f_i$ has rank less than or
equal to $r_1-1$, and similarly for the $g_i$. Then we require simply that
all of these subschemes be empty.
\end{Ilist}
\end{defn}

The main theorem of \cite{os8} on the linked Grassmannian is the following:

\begin{thm}\label{lg-main} \cite[Lem. A.3, Lem. A.12, Thm. A.15]{os8} 
$\cLG$ is representable by a
scheme $LG$; this scheme is naturally a closed subscheme of the obvious
product $G_1 \times \cdots \times G_n$ of Grassmannian schemes over $S$,
which is smooth of relative dimension $nr(d-r)$. Each component of $LG$ has
codimension $(n-1)r(d-r)$ in the product, and maps
surjectively to $S$; in fact, the smooth points of $LG$ are dense in every
component of every fiber. If $s$ is non-zero, then $LG$ is also irreducible.
\end{thm}

\section{The main calculation}

We now give an explicit local computation of the structure of a linked
Grassmannian of length $2$.

\begin{notn} Let $M_{d,S}(2,s)$ be the affine scheme over $S$ parametrizing
pairs of $d\times d$ matrices $M,N$ with coefficients in $\cO_S$ and with
$MN=NM=s$.
\end{notn}

\begin{thm}\label{calc} Fix vector bundles $\cE_1,\cE_2$ of rank $d$ on $S$, 
maps $f:\cE_1 \to \cE_2$ and $g:\cE_2 \to \cE_1$, $r \in \NN$ and 
$s \in \cO_S$ satisfying the conditions for
a linked Grassmannian of length $2$. Then given any field $k$ and a 
$k$-valued point $x=(V_1,V_2) \in LG:=LG(r,\{\cE_1,\cE_2\},\{f,g\})$, 
there exists a Zariski neighborhood $U$ of $x$ in $LG$ and an isomorphism
$$U \risom \AA^m_{S'} \times M_{\ell,S'}(2,s)$$
sending $x$ to the origin, where $S' \subseteq S$ is the image of $U$,
$$\ell=\begin{cases}r-\dim (g(V_2)) - \dim(f(V_1)): &s=0 \text{ at }x\\
0: &s \neq 0 \text{ at }x\end{cases},$$
and $m=r(d-r)-\ell^2$.
\end{thm}

\begin{proof}
The question is local on the base, so we may assume that $S = \Spec R$
and that $\cE_1$ and $\cE_2$ are trivial bundles. Furthermore, if 
$s \neq 0$ at $x$, then $LG$ is locally isomorphic simply to the 
Grassmannian $G(r,d)$, so in a neighborhood of $x$ we find that $LG$
is isomorphic to $\AA^{r(d-r)}_S$, as asserted. So we may suppose that
$s=0$ at $x$.

Let $\ox$ be the image of $x$ under the structure morphism from $LG$
to $S$.  Then the spaces $(V_1,V_2)$ determining $x$ are naturally 
subspaces of $\ox^*\cE_1$ and $\ox^*\cE_2$, and must satisfy
$f(V_1) \subseteq V_2$
and $g(V_2) \subseteq V_1$.  Let $c$ be the dimension of
the kernel of $f$ on $\ox^*\cE_1$, and let
$d_1$ and $d_2$ be the dimensions of $f(V_1)$ and $g(V_2)$, respectively.

We make the following choices:

\begin{itemize}
\item $\{\alpha^1_1, \dots, \alpha^1_{d_1}\}$ in $V_1$ so that 
$\{f(\alpha^1_1), \dots, f(\alpha^1_{d_1})\}$ is a basis for $f(V_1)$;
\item $\{\beta^1_1, \dots, \beta^1_{d_2}\}$ in $V_2$ so that
$\{g(\beta^1_{1}), \dots, g(\beta^1_{d_2})\}$ is a basis for $g(V_1)$; 
\item $\{\alpha^2_{1}, \dots, \alpha^2_{r - d_1 - d_2}\}$ in $\ox^*\cE_1$ so
that 
$\{f(\alpha^2_{1}), \dots, f(\alpha^2_{r - d_1 - d_2})\}$ 
extends the set
$\{f(\alpha^1_1), \dots, f(\alpha^1_{d_1})\}$ 
to a basis for $\ker g \cap V_2$;
\item $\{\beta^2_{1}, \dots, \beta^2_{r - d_1 - d_2}\}$ in $\ox^*\cE_2$ so
that 
$\{g(\beta^2_{1}), \dots, g(\beta^2_{r - d_1 - d_2})\}$ 
extends the set
$\{g(\beta^1_{1}), \dots, g(\beta^1_{d_2})\}$ 
to a basis for $\ker f \cap V_1$,
\item $\{\alpha^3_{1}, \dots, \alpha^3_{d - r - c + d_2}\}$ in $\ox^*
\cE_1$ so that
$\{f(\alpha^3_{1}), \dots, f(\alpha^3_{d-r-c-d_2})\}$ 
extends the set
$\{f(\alpha^1_1), \dots, f(\alpha^1_{d_1}), f(\alpha^2_{1}), \dots,
f(\alpha^2_{r - d_1 - d_2})\}$ 
to a basis for $\ker g$, and
\item $\{\beta^3_{1}, \dots, \beta^3_{d_1+c-r}\}$ in $\ox^* \cE_2$ so that
$\{g(\beta^3_{1}), \dots, g(\beta^3_{d_1+c-r})\}$ extends the set
$\{g(\beta^1_{1}), \dots, g(\beta^1_{d_2}),g(\beta^2_{1}),\dots,
g(\beta^2_{r-d_1-d_2})\}$ 
to a basis for $\ker f$.
\end{itemize}

We lift each $\alpha^i_j$ defined above to a section $\talpha^i_j$ of $\cE_1$,
and each $\beta^i_j$ to a section $\tbeta^i_j$ of $\cE_2$.  These lifts
allow us to define nice bases for these bundles in a neighborhood of $\ox$.
Specifically, for $i$ between $1$ and $d$, 
we define $e^1_i$ and $e^2_i$ as follows:

\begin{itemize}
\item $e^1_1,\dots,e^1_{d_1}:=\talpha^1_1,\dots,\talpha^1_{d_1}$;
\item $e^1_{d_1+1},\dots,e^1_{r-d_2}:=g(\tbeta^2_{1}),\dots,
g(\tbeta^2_{r-d_1-d_2})$;
\item $e^1_{r-d_2+1},\dots,e^1_{r}:=g(\tbeta^1_1),\dots,
g(\tbeta^1_{d_2})$;
\item $e^1_{r+1},\dots,e^1_{2r-d_1-d_2}:=\talpha^2_1,\dots,
\talpha^2_{r-d_1-d_2}$;
\item $e^1_{2r-d_1-d_2+1},\dots,e^1_{d+r-d_1-c}:=\talpha^3_1,\dots,
\talpha^3_{d-r-c+d_2}$;
\item $e^1_{d+r-d_1-c+1},\dots,e^1_d:=g(\tbeta^3_1),\dots,
g(\tbeta^3_{d_1+c-r})$;
\item $e^2_1,\dots,e^2_{d_1}:=f(\talpha^1_1),\dots,f(\talpha^1_{d_1})$;
\item $e^2_{d_1+1},\dots,e^2_{r-d_2}:=f(\talpha^2_1),\dots,
f(\talpha^2_{r-d_1-d_2})$;
\item $e^2_{r-d_2+1},\dots,e^2_{r}:=\tbeta^1_{1},\dots,\tbeta^1_{d_2}$;
\item $e^2_{r+1},\dots,e^2_{2r-d_1-d_2}:=\tbeta^2_{1},\dots,
\tbeta^2_{r-d_1-d_2}$;
\item $e^2_{2r-d_1-d_2+1},\dots,e^2_{d+r-d_1-c}:=f(\talpha^3_1),\dots,
f(\talpha^3_{d-r-c+d_2})$;
\item $e^2_{d+r-d_1-c},\dots,e^2_d:=\tbeta^3_{1},\dots,\tbeta^3_{d_1+c-r}$.
\end{itemize}

Then $e^1_1, \dots, e^1_d$ and $e^2_1, \dots, e^2_d$ restrict to bases
for $\ox^*\cE_1$ and $\ox^*\cE_2$, respectively. Shrinking $S$ to a 
smaller neighborhood of $\ox$, we can thus assume that 
$e^1_1, \dots, e^1_d$ is a basis for $\cE_1$
and $e^2_1, \dots, e^2_d$ is a basis for $\cE_2$.

With respect to these bases, the maps $f$ and $g$ have block forms
$$f = \begin{pmatrix}
\Id_{d_1} & 0 & 0 & 0 & 0 & 0\\
0 & 0 & 0 & \Id_{r - d_1 - d_2} & 0 & 0\\
0 & 0 & s\Id_{d_2} & 0 & 0 & 0\\
0 & s\Id_{r - d_1 - d_2} & 0 & 0 & 0 & 0\\
0 & 0 & 0 & 0 & \Id_{d - r - c + d_2} & 0\\
0 & 0 & 0 & 0 & 0 & s\Id_{d_1 + c - r}
\end{pmatrix}$$
$$g = \begin{pmatrix}
s\Id_{d_1} & 0 & 0 & 0 & 0 & 0\\
0 & 0 & 0 & \Id_{r - d_1 - d_2} & 0 & 0\\
0 & 0 & \Id_{d_2} & 0 & 0 & 0\\
0 & s\Id_{r - d_1 - d_2} & 0 & 0 & 0 & 0\\
0 & 0 & 0 & 0 & s\Id_{d - r - c + d_2} & 0\\
0 & 0 & 0 & 0 & 0 & \Id_{d_1 + c - r}
\end{pmatrix}.$$

If we consider $LG$ as a subscheme of $G(r,d)_R \times_R G(r,d)_R$ in
the obvious way, $x$ has an affine neighborhood in this product given
by a pair of $d-r$ by $r$ matrices $A^1$, $A^2$, where this pair
corresponds to the pair of subspaces $(V_1^{\prime}, V_2^{\prime})$
satisfying
$$V_1^{\prime} = \colspan \begin{pmatrix} \Id_r \\ A^1 \end{pmatrix}$$
$$V_2^{\prime} = \colspan \begin{pmatrix} \Id_r \\ A^2 \end{pmatrix}.$$

We wish to find equations for the intersection of this neighborhood with
$LG$.  A pair $(A^1, A^2)$ will lie in $LG$ if and only if 
$f(V_1^{\prime})$
is contained in $V_2^{\prime}$ and $g(V_2^{\prime})$ is contained in
$V_2^{\prime}$.  To write these conditions explicitly, we break the
matrices defining $V_1^{\prime}$ and $V_2^{\prime}$ into block form
as follows:

$$V_1^{\prime} = \colspan
\begin{pmatrix} 
\Id_{d_1} & 0 & 0\\
0 & \Id_{r - d_1 - d_2} & 0\\
0 & 0 & \Id_{d_2}\\
A^1_{11} & A^1_{12} & A^1_{13}\\
A^1_{21} & A^1_{22} & A^1_{23}\\
A^1_{31} & A^1_{32} & A^1_{33}
\end{pmatrix}
$$
and similarly for $A^2$.  Here the horizontal divisions of $A$ are chosen
so that the block division above corresponds to the block division for the
matrices for $f$ and $g$.

Now we have:
$$f(V_1^{\prime}) = \colspan
\begin{pmatrix}
\Id_{d_1} & 0 & 0\\
A^1_{11} & A^1_{12} & A^1_{13}\\
0 & 0 & s\Id_{d_2}\\
0 & s\Id_{r - d_1 - d_2} & 0\\
A^1_{21} & A^1_{22} & A^1_{23}\\
sA^1_{31} & sA^1_{32} & sA^1_{33}
\end{pmatrix}
$$
and it is easy to see that this is contained in $V_2^{\prime}$ if and only
if the matrix identity:
$$
\begin{pmatrix}
A^2_{11} & A^2_{12} & A^2_{13}\\
A^2_{21} & A^2_{22} & A^2_{23}\\
A^2_{31} & A^2_{32} & A^2_{33}
\end{pmatrix}
\begin{pmatrix}
\Id_{d_1} & 0 & 0\\
A^1_{11} & A^1_{12} & A^1_{13}\\
0 & 0 & s\Id_{d_2}
\end{pmatrix}
=
\begin{pmatrix}
0 & s\Id_{r - d_1 - d_2} & 0\\
A^1_{21} & A^1_{22} & A^1_{23}\\
sA^1_{31} & sA^1_{32} & sA^1_{33}
\end{pmatrix} 
$$
is satisfied.  Similarly, $g(V_2^{\prime})$ is contained in $V_1^{\prime}$
if and only if we have
$$
\begin{pmatrix}
A^1_{11} & A^1_{12} & A^1_{13}\\
A^1_{21} & A^1_{22} & A^1_{23}\\
A^1_{31} & A^1_{32} & A^1_{33}
\end{pmatrix}
\begin{pmatrix}
s\Id_{d_1} & 0 & 0\\
A^2_{11} & A^2_{12} & A^2_{13}\\
0 & 0 & \Id_{d_2}
\end{pmatrix}
=
\begin{pmatrix}
0 & s\Id_{r - d_1 - d_2} & 0\\
sA^2_{21} & sA^2_{22} & sA^2_{23}\\
A^2_{31} & A^2_{32} & A^2_{33}
\end{pmatrix} 
$$
If we solve for the matrix variables $A^1_{13}, A^1_{23}, A^2_{31},
A^2_{32}, A^2_{33}, A^2_{11}, A^2_{21},$ and $A^1_{22}$ in terms
of the remaining variables, all of the
above relations disappear except for the relations
$$A^1_{12}A^2_{12} = A^2_{12}A^1_{12} = s\Id_{r - d_1 - d_2},$$
and the result follows.
\end{proof}

\begin{rem} Although the proof of the theorem is carried out by choosing
bases, without resorting to such vulgarities we can at least see how one 
obtains from a $(V_1',V_2')$ in a
neighborhood of $(V_1,V_2)$ a pair of maps $\bar{f},\bar{g}$, between
$(r-d_1-d_2)$-dimensional vector spaces, such that $\bar{f} \circ \bar{g}=
\bar{g} \circ \bar{f} = s$. For convenience, we work over $\Spec k$, and
hence suppose $s=0$. As in the identification of a neighborhood of
a point in a Grassmannian with affine space, we do need to choose
complementary subspaces: thus, if our point $x$ corresponds to $(V_1,V_2)$,
with $V_i \subseteq \ox^* \cE_i$, we choose $W_i$ complementary to $V_i$
in $\ox^* \cE_i$ for $i=1,2$. We therefore obtain an affine neighborhood
of $(V_1,V_2) \in G(r,d) \times G(r,d)$ parametrized by 
$\Hom(V_1,W_1) \times \Hom(V_2,W_2)$.

Given $(V_1',V_2')$ corresponding to 
$(\psi_1,\psi_2) \in \Hom(V_1,W_1) \times \Hom(V_2,W_2)$, we will
construct $\bar{f}:W_1/(f^{-1}(W_2)\cap W_1) \to W_2/(g^{-1}(W_1)\cap W_2)$
and $\bar{g}:W_2/(g^{-1}(W_1)\cap W_2) \to W_1/(f^{-1}(W_2)\cap W_1)$
such that if $(V_1',V_2')$ are linked under $f,g$, then 
$\bar{f} \circ \bar{g}=\bar{g}\circ \bar{f}=s=0$. Note that these spaces
both have dimension $r-d_1-d_2$. The key point is to observe
that we have a natural isomorphism
$$\phi_1: f^{-1}(V_2) \cap W_1/(\ker f \cap W_1) 
\risom W_1/(f^{-1}(W_2)\cap W_1),$$
and similarly 
$$\phi_2: g^{-1}(V_1) \cap W_2/(\ker g \cap W_2) 
\risom W_2/(g^{-1}(W_1)\cap W_2),$$
induced by the inclusion maps. We can then define 
$\bar{f}= \psi_2 \circ f \circ \phi_1^{-1}$, and
$\bar{g}= \psi_1 \circ g \circ \phi_2^{-1}$, modding out appropriately in
the target space. 

It remains to check that if $(V_1',V_2')$ are linked under $f,g$, then
$\bar{f} \circ \bar{g}=\bar{g}\circ \bar{f}=s=0$. If we denote by
$\pr_V$ the projection to a subspace $V$, we observe that 
$f \circ \phi_1^{-1} \circ \psi_1 = \pr_{V_2} \circ f \circ \psi_1$,
so 
$\pr_{W_2} \circ f \circ \psi_1 = 
f \circ \psi_1-f \circ \phi_1^{-1} \circ \psi_1$.
Now, if $f(V_1') \subseteq V_2'$, we have by definition that for all 
$v \in V_1$,
$f(v)+f\psi_1(v)=u+\psi_2(u)$ for some $u \in V_2$; more specifically,
we must have $u=\pr_{V_2} (f(v)+f \psi_1(v))$, so we have 
$$f(v)+f\psi_1(v)=
\pr_{V_2}(f(v)+f \psi_1(v))+\psi_2 \pr_{V_2} (f(v)+f \psi_1(v)).$$ 
We have
that $f(v) \in V_2$, and $\psi_2$ takes values in $W_2$; if we 
apply $\pr_{W_2}$ to both sides of the equation, and substitute, we find 
$$f (\psi_1(v))-f (\phi_1^{-1} (\psi_1(v)))=
\psi_2(f(v))+\psi_2(f(\phi_1^{-1}(\psi_1(v))))$$ 
for all $v \in V_1$.
Finally, we note that $g \circ \phi_2^{-1}$ takes values in $V_1$,
and that $\bar{f}\circ\bar{g}$ takes values modulo $g^{-1}(W_1) \supseteq
\ker g = \im f$, so we see that 
$\psi_2 \circ f \circ \phi_1^{-1} \circ \psi_1 \circ g \circ \phi_2^{-1} =
\psi_2 \circ f \circ g \circ \phi_2^{-1} +
\psi_2 \circ f \circ \phi_1^{-1} \circ \psi_1 \circ g \circ \phi_2^{-1}$
takes values in $\im f$, which
shows that $\bar{f} \circ \bar{g}=0$, as desired. Similarly, if 
$g(V_2') \subseteq V_1'$, we check that $\bar{g} \circ \bar{f}=0$.
\end{rem}

\section{Conclusion of the main theorems}

Zhang \cite{zh1} has shown, using the theory of Hodge algebras of de Concini,
Eisenbud and Procesi \cite{d-e-p},
that spaces such as those arising in Theorem \ref{calc} 
are Cohen-Macaulay. We are thus able to show the following.

\begin{thm}\label{lg-flat} Let $S$ be integral and Cohen-Macaulay, 
and let $LG:=LG(r,\{\cE_i\}_i,\{f_i,g_i\}_i)$
be any linked Grassmannian scheme over $S$. Then $LG$ is flat over $S$,
and reduced and Cohen-Macaulay, with reduced fibers.
\end{thm}

\begin{proof} We first show that $LG$ is Cohen-Macaulay whenever 
$S=\Spec k$.
We prove this by induction on $n$, starting with the case $n=2$.
In this case, by Theorem 
\ref{calc}, we can conclude that $LG$ is Cohen-Macaulay if we know that 
$M_{\ell,S}(2,s)$ is Cohen-Macaulay. If $s \neq 0$, we have $\ell=0$, and 
$M_{\ell,S}(2,s)$ is 
just a point. If $s=0$, we have the space of pairs of $d\times d$ matrices
$M,N$ with coefficients in $k$ and with $MN=NM=0$. By \cite[Thm. 22]{zh1},
this space is Cohen-Macaulay. 

We now induct on $n$. 
We let $LG$ be a linked Grassmannian of length $n$, and denote by 
$LG_{n-1}$ (respectively, $LG'$) the truncated linked Grassmannians of 
length $n-1$ (respectively, $2$) obtained from $LG$ by forgetting 
$\cE_n,f_n,g_n$ (respectively, $\cE_i,f_i,g_i$ for all $i<n-1$).
Then $LG_{n-1}$ naturally lies in the product of Grassmannians
$G_{n-1}:=G(r,\cE_1)\times \dots \times G(r,\cE_{n-1})$, and
$LG'$ in $G':=G(r,\cE_{n-1}) \times G(r,\cE_n)$.
Let 
$$\Delta:G:=G(r,\cE_1)\times\dots \times G(r,\cE_n) \to 
G_{n-1} \times G'$$ 
be the map which doubles the $(n-1)$st coordinate; it is then clear by
considering the functors of points that we have that
$LG \subseteq G$ is the fiber product of $\Delta$ with the inclusion map
$LG_{n-1} \times LG' \to G_{n-1} \times G'$.

Now, the image of $\Delta$ is a locally complete intersection of codimension
$r(d-r)$ in $G_{n-1}\times G'$, so we conclude that $LG$ is cut out
by $r(d-r)$ equations inside of $LG_{n-1} \times LG'$, and since this is
also the codimension of $LG$, we find that the induction hypothesis that
$LG_{n-1}$ and $LG'$ are Cohen-Macaulay then implies that $LG$ is 
Cohen-Macaulay as well. This shows that $LG$ is always Cohen-Macaulay over
$\Spec k$. 

Next, still assuming that $S=\Spec k$, by Theorem \ref{lg-main}, we have 
that $LG$ is smooth and hence reduced on an open dense subset, and since 
$LG$ is Cohen-Macaulay, it cannot have imbedded components, and we conclude 
that $LG$ is also reduced as long as the base is a point. 

Finally, suppose $S$ is an arbitrary integral, Cohen-Macaulay scheme. The
base change of $LG$ to any integral, Cohen-Macaulay scheme is another linked 
Grassmannian, so since we have reducedness when $S=\Spec k$, we have that 
condition (i) of the following lemma is satisfied for any base.
Similarly, by Theorem \ref{lg-main} every topological component of $LG$ 
surjects onto the base after base change to $\Spec A$ with $A$ a DVR,
and condition (ii) of the following lemma is satisfied. Thus, the lemma 
implies that any linked Grassmannian is reduced and flat. Finally, we can 
conclude Cohen-Macaulayness for arbitrary $S$ by \cite[Cor., p. 181]{ma1}. 
\end{proof}

\begin{lem} Let $f:X \to S$ be a morphism locally of finite presentation
with $S$ reduced and Noetherian, and suppose that:
\begin{ilist}
\itm the fibers of $f$ are geometrically reduced;
\itm after arbitrary base change to $\Spec A$, with $A$ a DVR, no 
(topological) component of $X \times _S \Spec A$ is supported over the 
special fiber.
\end{ilist}
Then $X$ is reduced, and $f$ is flat.
\end{lem}

\begin{proof} By \cite[Thm. 11.8.1]{ega43},
we can prove that $f$ is flat by checking that it is flat after arbitrary
base change under morphisms $\Spec A \to S$, where $A$ is a DVR. Since
$f$ has geometrically reduced fibers, the fibers remain reduced under base
change, and assumption (ii) allows us to conclude that 
$(X \times _S \Spec A)_{\red}$ is flat over $\Spec A$ by 
\cite[Prop. 14.3.8]{ega43}.

It then follows by \cite[Lem. 6.13]{os8} that $X \times _S \Spec_A$ is
reduced, hence flat, giving us that $f$ is flat. The flatness of $f$
then also implies by \cite[Prop. 15.3.1]{ega43} that $X$ is reduced,
as desired.
\end{proof}

Finally, we use the previous theorem to conclude a similar flatness result 
for the limit linear series schemes introduced in \cite{os8}; see that
paper for the appropriate definitions. 

\begin{thm}\label{main-flat} Fix integers $r,d,n$, and ramification sequences
$\alpha^1,\dots,\alpha^n$, and let $X_0$ be a 
curve of compact type over $\Spec k$ having no more than two components,
with marked points $P_1,\dots,P_n$. Suppose that 
the space $G^r_d(X_0):=G^r_d(X_0/k; \{(P_i,\alpha^i)\}_i)$
of limit linear series 
on $X_0$ of degree $d$ and dimension $r$, with prescribed ramification at the
$P_i$, has exactly the expected dimension 
$\rho:=(r+1)(d-r)-rg-\sum_{i,j} \alpha^i_j$. 

Then $G^r_d(X_0)$ is Cohen-Macaulay, and if we have a 
smoothing family of curves $X/B, \{\tilde{P}_i\}_i$ (see \cite[Def. 3.1]{os8}) 
where $X_0, \{P_i\}_i$ occurs as a fiber over $b_0 \in B$, then 
$G^r_d(X):=G^r_d(X/B;\{(\tilde{P}_i,\alpha^i)\}_i)$, the relative space of 
limit linear series on $X/B$, is Cohen-Macaulay and flat over $B$ at every 
point over $b_0$. 
\end{thm}

We remark that the expected dimension hypothesis is satisfied for a 
general curve with marked points in characteristic 0; see \cite{os13}.

\begin{proof} Indeed, in the proof of \cite[Thm. 5.3]{os8} the 
space $G^r_d(X)$ of limit linear series on a
smoothing family $X/B$ (where we may allow $B=\Spec k$, $X=X_0$), is
constructed inside a product of Grassmannians by intersecting a linked
Grassmannian with a finite collection of Schubert cycles, and this has
dimension $\rho$ in the fiber over $b_0$ if and only if this intersection
has maximal codimension over $b_0$. Since we now know that linked
Grassmannians are Cohen-Macaulay, at any point over $b_0$ we have a 
proper-dimensional intersection of Cohen-Macaulay schemes inside a regular 
one, and this must be Cohen-Macaulay by the lemma below.
Since $B$ is regular as part of the hypothesis of a smoothing
family, we also obtain that $G^r_d(X)$ is flat over $B$ (at least, at any
point over $b_0$) by \cite[Thm. 23.1]{ma1}.
\end{proof}

We recall the following well-known lemma, and include the proof for lack
of a suitable reference:

\begin{lem}\label{cm-int} Let $X$ be regular, and $Y_1,\dots,Y_n$ closed
subschemes of $X$, each Cohen-Macaulay, of codimensions $c_1,\dots,c_n$.
Let $Z=Y_1 \cap \dots \cap Y_n$, and suppose that $Z$ has codimension
$\sum_i c_i$ in $X$. Then $Z$ is also Cohen-Macaulay.
\end{lem}

\begin{proof} By induction, it is clearly enough to consider the case
that $n=2$. Also, the question is purely local, so we are reduced to
the following: let $R$ be a regular local ring, and $I,J$ ideals of 
codimension $d,e$ such that $R/I$ and $R/J$ are Cohen-Macaulay, and $I+J$
has codimension $d+e$. Then we want to see that $R/(I+J)$ is Cohen-Macaulay.
By \cite[Cor. 19.15]{ei1}, the Cohen-Macaulayness of the rings 
$R/I,R/J,R/(I+J)$ is equivalent to having free resolutions over $R$ of
lengths $d,e,d+e$ respectively. Thus, we may assume we have free resolutions
$\cF^{\bullet},\cG^{\bullet}$ of $R/I$ and $R/J$ having lengths $d,e$
respectively. We claim that that the tensor complex $\cF^{\bullet} \otimes
\cG^{\bullet}$, which has length $d+e$, is still exact, thus giving a
free resolution of $R/(I+J)$ of length $d+e$, and showing that $R/(I+J)$
is Cohen-Macaulay.

By \cite[Thm. 20.9]{ei1}, it is enough to check that 
$\cF^{\bullet} \otimes \cG^{\bullet}$
is exact after localizing at an arbitrary prime $\fp$ of codimension 
strictly less than $d+e$: indeed, condition (1) of {\it loc.\ cit.} is 
automatically satisfied for the tensor product of two free resolutions, and
to verify (2), if we had (in the notation of {\it loc.\ cit.}) 
$\depth I(\vp_k)<k$ for some $k \leq d+e$, 
we could choose $\fp$ to be a minimal prime of $I(\vp_k)$ having codimension
strictly less than $d+e$, and (2) would still be violated after localizing at
$\fp$. 

Now, let $\fp$ be arbitrary of codimension less than $d+e$; it cannot 
contain $I+J$, and thus does
not contain both $I$ and $J$; without loss of generality, say that $I$ is
not contained in $\fp$. But if we localize $\cF^{\bullet}$ at $\fp$,
$R/I$ itself becomes $0$, so the sequence becomes split exact by 
\cite[Lem. 20.1]{ei1}, and in this case, tensoring
with $(\cG^{\bullet})_{\fp}$ necessarily preserves exactness.
\end{proof}

\bibliographystyle{hamsplain}
\bibliography{hgen}

\end{document}